\documentclass[12]{amsart}
\usepackage{latexsym}
\usepackage{amsthm}
\usepackage{amssymb}
\usepackage{amsfonts, amsmath}

\newtheorem{lm}{Lemma}

\newtheorem{Thm}{Theorem}

\title{On $q$-complete and $q$-concave with corners complex manifolds.}
\author{}
\begin{document}
\maketitle
\begin{center}
\em{Youssef Alaoui}\\
\em{y.alaoui@iav.ac.ma}\\
\end{center}

\noindent
\em{Department of Mathematics,
Hassan II Institute of Agronomy}\\ 
\em{and Veterinary Sciences,
Madinat Al Irfane, BP 6202, Rabat, 10101, Morocco,}\\
\date{}
\linespread{1.3}
\date{}
\begin{abstract}
It is proved that if there exists a positive and continuous function $f$ on an
$n$-dimensional complex manifold $X$, $q$-convex with corners outside a compact set $K\subset X$ and which exhausts $X$ from below, then $dim_{\mathbb{C}}H^{p}(X,{\mathcal{F}})<+\infty$ for any coherent analytic sheaf ${\mathcal{F}}$ on $X$ if $p<n-q$.\\
\hspace*{.1in}It is known from the theory of Andreotti and Grauert that if a complex space
$X$ is $q$-complete, then $X$ is cohomoloogically $q$-complete.
Until now it is not known in general if these two conditions are equivalent.\\
The aim of section $4$ of this article is to provide a counterexample to the conjecture posed by Andreotti and Grauert ~\cite{ref2} to show that a cohomologically $q$-complete space is not necessarily $q$-complete.\\
\hspace*{.1in}In section $5$ of this article, we will prove
that there exist for each pair of integers $(n,q)$ with $2\leq q\leq n-1$
a $q$-complete with corners open subset $D$ of $\mathbb{P}^{n}$ and $\mathcal{F}\in coh(\mathbb{P}^{n})$ such that $D$ is not cohomologically
$\hat{q}$-complete with respect to ${\mathcal{F}}$. Here $\hat{q}=n-[\frac{n-1}{q}]$, where $[x]$ denotes the integral part of $x$.
\end{abstract}

\section{Introduction}
Finiteness and Vanishing theorems of Andreotti and Grauert ~\cite{ref2} play a very
impotant role in the theory of Complex Analytic Geometry. These theorems follow from the existence of smooth $q$-convex and $q$-concave exhaustion functions.\\
In many examples, the natural exhaustion function $f$ is not smooth but only
locally the maximum of finitely many $q$-convex functions $f=max(f_{1},\cdots,f_{s})$.\\
\hspace*{.1in} In ~\cite{ref5}, Diederch and Fornaess have proved that every $q$-convex with corners function on a complex manifold of dimension $n$ can be approximated in $C^{0}$
topology by $\tilde{q}$-convex functions on $X$, where $\tilde{q}=n-[\frac{n}{q}]+1$.
They moreover showed by means of a counter-example that the number $\tilde{q}$ obtained is optimal.\\
\hspace*{.1in}It was shown by Andreotti and Grauert ~\cite{ref3} that if $X$
is a $q$-concave complex space, then for any ${\mathcal{F}}\in coh(X)$,
$dim_{\mathbb{C}}H^{p}(X,{\mathcal{F}})<+\infty$ if $p<prof({\mathcal{F}})-q$.\\
\hspace*{.1in}In section $3$ of this paper, we prove an extension of this result
for families of finitely dimensional $q$-concave with corners complex manifolds. (For the definitions, see below.)\\
\hspace*{.1in}
In 1962, A. Andreotti and H. Grauert ~\cite{ref2} showed finiteness and vanishing theorems for
cohomology groups of analytic spaces under geometric conditions of q-convexity. Since then
the question whether the reciprocal statements of these theorems are true have been subject to
extensive studies, where for $q > 1$ more specific assumptions have been added. For example, it
is known from the theory of Andreotti–Grauert that a q-complete complex space is always
cohomologically q-complete, but it is not known if these two conditions are equivalent except
when $X$ is a Stein manifold, $\Omega \subset X$ is cohomologically $q$-complete with respect to $\mathcal{O}_{\Omega}$ and $\Omega$
has a smooth boundary ~\cite{ref3}.\\

\hspace*{.1in}In the present article, Section 4 is devoted to establishing a counter-example to the Andreotti and Grauert conjecture. This construction is explicit and constructive in nature. Specifically, we show the existence of a connected closed submanifold $A\subset\mathbb{P}^5$ of codimension $3$ such that $\mathbb{P}^5\setminus A$ is cohomologically 3-complete but not 3-complete.\\
\hspace*{.1in}In ~\cite{ref8}, Matsumoto has shown the vanishing theorems
for an intersection of a finite number of $q$-complete domains in a complex manifold
of dimension $n$. She has proved that if $D_{1},\cdots, D_{t}$
are $q$-complete open subsets of a complex manifold $M$ of dimension $n$ and,
if ${\mathcal{F}}$ is a coherent analytic sheaf on $M$ such that $H^{n}(M,{\mathcal{F}})=0$,
then
$$H^{p}(D_{1}\cap\cdots\cap D_{t},{\mathcal{F}})=0 \ \ \text{for \ all} \ \ p\geq \hat{q},$$
where
\[\hat{q}=n-[\frac{n-1}{q}]=\left\{\begin{array}{cc}
\tilde{q} \ \ if \ \  q\mid n\\
\tilde{q}-1 \ \ if \ \ q \nmid n
\end{array}
\right.\]
Here $\tilde{q}=n-[\frac{n}{q}]+1$ and $\hat{q}=n-[\frac{n-1}{q}]$, where $[x]$
denotes the integral part of $x$.\\
But it is not known if the same result follows if $D$ is an arbitrary $q$-complete
with corners open subset of $M$.\\
\hspace*{.1in}In section $5$ of this paper, we will prove by means of a counterexample
that there exist for each pair of integers $(n,q)$ with $2\leq q\leq n-1$
a $q$-complete with corners open subset $D$ of $\mathbb{P}^{n}$ and
$\mathcal{F}\in coh(\mathbb{P}^{n})$ such that $H^{\hat{q}}(D,{\mathcal{F}})\neq 0$ .
\section{Preliminaries}
We start by recalling some definitions and results concerning $q$-convexity.\\
\hspace*{.1in}Let $X$ be a complex manifold. Then it is known that
a function $\phi\in C^{\infty}(X)$ is $q$-convex if for every
point $z\in X$, the Levi form
$L_{z}(\phi; z)$ has at most $q-1$ negative or zero eingenvalues
on each tangent space $T_{z}\Omega$, $z\in X$.\\
\hspace*{.1in}We say that $X$ is $q$-complete if there exists
a $q$-convex function $\phi\in C^{\infty}(X,\mathbb{R})$
which is exhaustive on $X$ i.e. $\{x\in X: \phi(x)<c\}$
is relatively compact in $X$ for any $c\in\mathbb{R}.$\\
\hspace*{.1in}The space $X$ is said to be cohomologically $q$-complete
if for every coherent analytic sheaf ${\mathcal{F}}$ on $X$ the cohomology groups $H^{r}(X, {\mathcal{F}})$ vanish for all $r\geq q$.\\
\hspace*{.1in}An open subset $D$ of $X$ is called $q$-Runge
if for every compact set $K\subset D$, there is a $q$-convex
exhaustion function $\phi\in C^{\infty}(X)$ such that
$$K\subset\{x\in X: \phi(x)<0\}\subset\subset D$$
This generalizes the classical notion of Runge pairs of Stein spaces.\\
It is shown in ~\cite{ref2} that if $D$ is $q$-Runge in $X$, then for
every ${\mathcal{F}}\in coh(X)$ the cohomology groups $H^{p}(D,
\mathcal{F})$ vanish for $p\geq q$ and, the restriction map
$$H^{p}(X, {\mathcal{F}})\longrightarrow H^{p}(D, {\mathcal{F}})$$
has dense image for all $p\geq q-1$.\\
\hspace*{.1in} A function $f : X\rightarrow \mathbb{R}$ is called $q$-convex
with corners , if $f$ is continuous and for each $x\in X$, there are a neighborhood
$U$ of $x$ in $X$ and $q$-convex functions $\phi_{1},\cdots, \phi_{r}$
on $U$ with $f|_{U}=max(\phi_{1},\cdots, \phi_{r})$.\\
We denote by $F_{q}(X)$ the set of the $q$-convex functions with corners on $X$.\\
\hspace*{.1in}A complex space $X$ will be called $q$-concave with corners if
there exists a continuous function $f : X\rightarrow ]0,+\infty[$ which is
$q$-convex with corners outside a compact set $K\subset X$ and such that $X_{c}=\{x\in X : f(x)>c\}\subset\subset X$ for each $c>0$.\\
\hspace*{.1in}The space $X$ is called $q$-complete with corners if there exists
a $q$-convex with corners exhaustion function $f\in F_{q}(X)$.
\section{$q$-concavity with corners}
\begin{lm}{Let $X$ be a complex manifold of dimension $n$, and let $\phi:
X\rightarrow \mathbb{R}$ be a smooth $q$-convex function $\phi$ on $X$. Let $\xi_{0}\in X$
and $X'_{c}=\{x\in X: \phi(x)>c\},$ where $c=\phi(\xi_{0})$.
Then for any coherent analytic sheaf ${\mathcal{F}}$
on $X$ the restriction map
$$H^{p}(X,{\mathcal{F}})\rightarrow H^{p}(X'_{c}, {\mathcal{F}})$$
is bijective if $p\leq n-q-1$,\\
\hspace*{.1in}injective if $p=n-q.$}
\end{lm}
Let $D$ be a domain in $\mathbb{C}^{n}$, $\xi\in D$,
and let $\phi\in C^{\infty}(D )$ be a q-convex
function. Then in order to prove lemma $1$ we shall need the following result due to Andreotti and Grauert ~\cite{ref3}.
\begin{Thm}{For any coherent analytic sheaf ${\mathcal{F}}$ on $D$
there exists a fundamental system of Stein
neighborhoods $U\subset D$ of $\xi$ such that if $Y=\{z\in D:
\phi(z)>\phi(\xi)\}$, then $H^{p}(Y\cap U, {\mathcal{F}})=0$ for
$0<p<n-q$ and $H^{0}(U,{\mathcal{F}})\rightarrow H^{0}(U\cap Y, {\mathcal{F}})$ is an isomorphism.}
\end{Thm}
\noindent
\begin{proof}
Let $V\subset\subset X$ be an open neighborhood of $\xi_{0}$ biholomorphic to a domain in $\mathbb{C}^{n}$.
Then there exists, by theorem $1$, a fundamental system of connected Stein neighborhoods $U\subset V$ of $\xi_{0}$ such that $H^{r}(U\cap X'_{c}, {\mathcal{F}})=0$ for $1\leq r<n-q$ and
$H^{0}(U, {\mathcal{F}})\rightarrow H^{0}(U\cap X'_{c}, {\mathcal{F}})$ is an isomorphism, or equivalently (See ~\cite{ref7} or
~\cite{ref1}), $\underline{H^{r}_{S}}({\mathcal{F}})=0$ for $r\leq
n-q,$ where
$\underline{H^{r}_{S}}({\mathcal{F}})$ is the cohomology sheaf
with support in $S=\{x\in X: \phi(x)\leq c\}$ and
coefficients in ${\mathcal{F}}.$ Furthermore, there exists a
spectral sequence
$$H^{p}_{S}(X, {\mathcal{F}})\Longleftarrow
E_{2}^{p,q}=H^{p}(X, \underline{H^{q}_{S}}({\mathcal{F}}))$$
Since $\underline{H^{p}_{S}}({\mathcal{F}})=0$ for $p\leq n-q,$ then for any $p\leq n-q,$
the cohomology groups $H^{p}_{S}(X, {\mathcal{F}})$ vanish and, the exact sequence of local
cohomology
\begin{center}
$\cdots\rightarrow H^{p}_{S}(X, {\mathcal{F}}) \rightarrow
H^{p}(X, {\mathcal{F}}) \rightarrow H^{p}(X'_{c}, {\mathcal{F}})
\rightarrow H^{p+1}_{S}(X, {\mathcal{F}})\rightarrow\cdots$
\end{center}
implies that $H^{p}(X, {\mathcal{F}}) \rightarrow H^{p}(X'_{c}, {\mathcal{F}})$
is bijective for any $c\in\mathbb{R}$ if $p\leq n-q-1$ and,
injective if $p=n-q$.
\end{proof}

\begin{lm}{Let $D$ be an open set in $\mathbb{C}^{n}$, $f\in F_{q}(D)$,
$1\leq q\leq n-1$. Then there exists for each point $\xi_{0}\in D$
a fundamental system of Stein neighborhoods $U$ of $\xi_{0}$ such that if $Y=\{z\in D : f(z)>f(\xi_{0})\}$, then for any coherent analytic sheaf ${\mathcal{F}}$ on $D$ we have :\\
(i) $H^{0}(U, \mathcal{F})\rightarrow H^{0}(U\cap Y, \mathcal{F})$ is bijective; \\
(ii) $H^{r}(U\cap Y, \mathcal{F})=0$ for $0<r<n-q$.}
\end{lm}
\begin{proof}
Let $U$ be a Stein neighborhood of $\xi_{0}$ in $D$ such that there exist
finitely many $q$-convex functions $\phi_{1},\cdots, \phi_{s} : U\rightarrow \mathbb{R}$
with $f|_{U}=max(\phi_{1},\cdots, \phi_{s})$.\\
By suitable choice of $U$, assertions (i) and (ii) are true when $s=1$, according to
theorem $1$. We, obviously, also may assume that the restriction
$f|_{U} : U\rightarrow \mathbb{R}$ is the maximum of two $q$-convex functions
$f|_{U}=max(\phi_{1},\phi_{2})$, which implies that $Y\cap U=Y_{1}\cup Y_{2},$
where $Y_{i}=\{z\in U : \phi_{i}(z)>f(\xi_{0})\}$ for $i=1, 2$.
If $q+1<n$, then by lemma $1$, we may choose $U$ so that
$H^{0}(U, \mathcal{F})\cong H^{0}(Y_{i}, \mathcal{F})$,
$H^{1}(Y_{i}, \mathcal{F})\cong H^{1}(U, \mathcal{F})=0$. Moreover, if for $j\geq 1$
the open set $Z_{j}=\{z\in Y_{1} : \phi_{2}(z)>f(\xi_{0})-\frac{1}{j}\}$ is not empty,
then by lemma 1, the restriction
$$H^{p}(Y_{1}, \mathcal{F})\rightarrow H^{p}(Z_{j}, \mathcal{F})$$
is bijective for $p=0, 1$. Therefore by  Mittag-Leffler theorem it follows that 
$$H^{p}(Y_{1}, \mathcal{F})=H^p\left(\lim _{\longleftarrow} Z_j, \mathcal{F}\right) \cong \lim _{\longleftarrow} H^p\left(Z_j, \mathcal{F}\right) \cong H^p(Y_{1} \cap Y_{2}, \mathcal{F}) \ \ \text{for} \ \ p=0, 1$$
This proves that
$H^{0}(Y_{i}, \mathcal{F})\cong H^{0}(Y_{1}\cap Y_{2}, \mathcal{F})$ 
and $H^{1}(Y_{i}, \mathcal{F})\cong H^{1}(Y_{1}\cap Y_{2}, \mathcal{F})=0$
It follows from the Mayer-Vietoris sequence for cohomology
\begin{center}
$0\rightarrow H^{0}(U\cap Y, \mathcal{F})\rightarrow H^{0}(Y_{1}, \mathcal{F})\oplus
H^{0}(Y_{2}, \mathcal{F})\rightarrow H^{0}(Y_{1}\cap Y_{2}, \mathcal{F})\rightarrow
H^{1}(U\cap Y, \mathcal{F})\rightarrow 0$
\end{center}
that $H^{1}(U\cap Y, \mathcal{F})=0$ and
$H^{0}(U, \mathcal{F})\cong H^{0}(U\cap Y, \mathcal{F})$.\\
\hspace*{.1in}Now if $2\leq r<n-q$, then by theorem $1$ we may take $U$
such that $H^{r-1}(Y_{i}, \mathcal{F})\cong H^{r}(Y_{i}, \mathcal{F})=0$ for $i=1, 2$
and, a proof similar to the one used previously shows that
$H^{r-1}(Y_{1}\cap Y_{2}, \mathcal{F})\cong H^{r-1}(Y_{1}, \mathcal{F})=0$, then the Mayer-Vietoris sequence for cohomology
\begin{center}
$\cdots\rightarrow H^{r-1}(Y_{1}, \mathcal{F})\oplus
H^{r-1}(Y_{2}, \mathcal{F})\rightarrow H^{r-1}(Y_{1}\cap Y_{2}, \mathcal{F})\rightarrow
H^{r}(U\cap Y, \mathcal{F})\rightarrow
H^{r}(Y_{1}, \mathcal{F})\oplus
H^{r}(Y_{2}, \mathcal{F})\rightarrow\cdots$
\end{center}
implies that  $H^{r}(U\cap Y, \mathcal{F})=0$.
\end{proof}
\begin{Thm}{Let $X$ be a $q$-concave with corners complex manifold of dimension $n$.
Then for any coherent analytic sheaf ${\mathcal{F}}$ on $X$ one has
$dim_{\mathbb{C}}H^{p}(X,{\mathcal{F}})<+\infty$ if $0\leq p< n-q.$}
\end{Thm}
\begin{proof}
The proof of theorem $2$ is similar to that of lemma $1$. In fact,
since $X$ is $q$-concave with corners, then there exists a continuous function
$f : X\rightarrow ]0,+\infty[$ which is $q$-convex with corners outside a compact set
$K\subset X$ and such that $X'_{c}=\{x\in X : f(x)>c\}\subset\subset X$ for every $c>0$.\\
\hspace*{.1in}Let $\xi_{0}\in X\setminus K$ be such that
$f(\xi_{0})=c$, and let $V\subset\subset X\setminus K$
be an open neighborhood of $\xi_{0}$ that can be identified with a domain of $\mathbb{C}^{n}$.
Then there exists, by lemma $2,$ a fundamental system of connected Stein neighborhoods $U\subset V$ of $\xi_{0}$ such that $H^{r}(U\cap X'_{c}, {\mathcal{F}})=0$ for $1\leq r<n-q$ and
$H^{0}(U, {\mathcal{F}})\rightarrow H^{0}(U\cap X'_{c},
{\mathcal{F}})$ is an isomorphism, which implies that if
$S=\{x\in X: \phi(x)\leq c\}$, then
the cohomology sheaf $\underline{H^{r}_{S}}({\mathcal{F}})=0$ for $r\leq n-q$.
Therefore for any $p\leq n-q,$
the cohomology groups $H^{p}_{S}(X, {\mathcal{F}})$ vanish and, the exact sequence of local cohomology
\begin{center}
$\cdots\rightarrow H^{p}_{S}(X, {\mathcal{F}}) \rightarrow
H^{p}(X, {\mathcal{F}}) \rightarrow H^{p}(X'_{c}, {\mathcal{F}})
\rightarrow H^{p+1}_{S}(X, {\mathcal{F}})\rightarrow\cdots$
\end{center}
yields that the map $H^{p}(X, {\mathcal{F}})\rightarrow H^{p}(X'_{c},
{\mathcal{F}})$ is bijective if $p\leq n-q-1$ 
and injective if $p=n-q$.
Since $X'_{c}\subset\subset X$, it follows from ~\cite{ref2} that
$dim_{\mathbb{C}}H^{p}(X, {\mathcal{F}})<+\infty$ if $p\leq n-q-1$.
\end{proof}
\section{A counterexample to the Andreotti-Grauert conjecture}
Let $A\subset \mathbb{P}^{n}$ be a closed submanifold of $codim(A)\leq q$.
Then by theorem $6$ of ~\cite{ref10} $\mathbb{P}^{n}\setminus A$ is $q$-complete
with corners. Consider the Veronese surface $A=\nu(\mathbb{P}^{2})$, where
$\nu : \mathbb{P}^{2}\rightarrow \mathbb{P}^{5}$ is the embedding given by
$$\nu([x: y: z])=[x^{2}: y^{2}: z^{2}: yz: xz: xy]$$
Then $\mathbb{P}^{5}\setminus A$ is $3$-complete with corners.
It was shown in ~\cite{ref6} that
$H^{5}(\mathbb{P}^{5}\setminus A, \mathbb{Z})=\mathbb{Z}/2\mathbb{Z}$.
By Morse theory it follows that $\mathbb{P}^{5}\setminus A$ is not $3$-complete.\\
\hspace*{.1in}By considering the resolution of the constant sheaf $\mathbb{C}$ given by :
$$0\rightarrow \mathbb{C}\rightarrow {\mathcal{O}}\stackrel{d}\rightarrow\Omega^{1}\stackrel{d}\rightarrow\Omega^{2}\rightarrow\cdots\stackrel{d}\rightarrow \Omega^{5}\rightarrow 0$$
where $\Omega^{i}$ denotes the sheaf of germs of holomorphic $p$-forms, and the fact
that $\mathbb{P}^{5}\setminus A$ is obviously cohomologically $3$-complete
with respect to the $\Omega^{i}$, we deduce that $\mathbb{P}^{5}\setminus A$ must
satisfies the condition $H^{p}(\mathbb{P}^{5}\setminus A, \mathbb{C})=0$
for all $p\geq 8$. It follows from a result due to Barth ~\cite{ref4} that
$\mathbb{P}^{5}\setminus A$ is cohomologically $3$-complete with respect to
coherent sheaves on $\mathbb{P}^{5}$. The mean purpose in this section is to prove that
$\mathbb{P}^{5}\setminus A$ is cohomologically $3$-complete; this gives a 
counterexample to the Andreotti-Grauert conjecture. (See ~\cite{ref2}).
\noindent
\begin{lm} Let $f \in F_3\left(\mathbb{P}^5 \backslash A\right)$ be a 3 -convex with corners exhaustion function on $\mathbb{P}^5 \backslash A$. There exists for each point $\xi_0 \in \mathbb{P}^5 \backslash A$ a Stein open neighborhood $U$ of $\xi_0$ such that if $Y=\left\{z \in \mathbb{P}^5 \backslash A: f(z)<f\left(\xi_0\right)\right\}$, then for any coherent analytic sheaf $\mathcal{F}$ on $U$ the cohomology group $H^p(U \cap Y, \mathcal{F})$ vanishes for all $p \geq 3$.
\end{lm}
\begin{proof}
Let $U \subset \subset \mathbb{P}^5 \backslash A$ be a Stein open neighborhood of $\xi_0$ such that there exist finitely many $q$-convex functions $\phi_1, \cdots, \phi_s: U \rightarrow \mathbb{R}$ with $\left.f\right|_U=\max \left(\phi_1, \cdots, \phi_s\right)$. 
Then $U \cap Y=Y_1 \cap \cdots \cap Y_s$, where $Y_i=\left\{z \in U: \phi_i(z)<f\left(\xi_0\right)\right\}, i=1, \cdots, s$, is $3$-complete and $3$-Runge  in $U$, because $U$ is Stein and $\phi_i$ is 3-convex on $U$. This implies that the restriction map

$$
H^p(U, \mathcal{F}) \rightarrow H^p\left(Y_i, \mathcal{F}\right)
$$
has a dense image if $p \geq 2$. Since by ~\cite{ref9} the canonical topologies on $H^p\left(Y_i, \mathcal{F}\right)$ are separated for all $p \geq 2$, then $H^p(Y_{i}, \mathcal{F})=0$ for all $i\in\{1,\cdots, s\}$ if $p \geq 2$. Therefore, if $s=2$, it follows from the mean theorem of
~\cite{ref8} that $H^p(Y_{1}\cap Y_{2}, \mathcal{F})=0$ for $p\geq 2\times 1+1=3$. Suppose now that $s\geq 3$ and for any $k$ with
$1\leq k\leq s-1$ the family $\{Y_{1},\cdots, Y_{s}\}$ satisfies the condition :
$$H^{p}(Y_{i_{1}}\cap\cdots\cap Y_{i_{k}}, {\mathcal{F}})=0$$
for all $p\geq 3$ and $i_{1}, i_{2},\cdots, i_{k}\in\{1, 2,\cdots, s\}$.
Then, by Proposition $1$ of ~\cite{ref8}, one obtains
$$H^{p}(Y_{1}\cap\cdots\cap Y_{s}, {\mathcal{F}})\cong H^{p+s-1}(Y_{1}\cup\cdots\cup Y_{s}, {\mathcal{F}})=0$$
for all $p\geq 3$, since $p+s-1\geq 5$. This completes the proof of lemma $3$.
\end{proof}
\begin{Thm}{The space $\mathbb{P}^{5}\setminus \nu(\mathbb{P}^{2})$, where
$\nu : \mathbb{P}^{2}\rightarrow \mathbb{P}^{5}$ is the Veronese embedding,
is not $3$-complete but for any coherent analytic sheaf
${\mathcal{F}}$ on $\mathbb{P}^{5}\setminus \nu(\mathbb{P}^{2})$
the cohomology group  $H^{p}(\mathbb{P}^{5}\setminus \nu(\mathbb{P}^{2}), {\mathcal{F}})$
vanishes for all $p\geq 3$.}
\end{Thm}
\begin{proof}
Let $f\in F_{3}(\mathbb{P}^{5}\setminus \nu(\mathbb{P}^{2}))$ be a $3$-convex with corners exhaustion function 
and denote by $X(\lambda)=\{z\in \mathbb{P}^{5}\setminus \nu(\mathbb{P}^{2}) : f(z)=\lambda\}$ for every $\lambda\in \mathbb{R}$.
We claim that for every pair of real numbers $\lambda<\mu$ we have:\\
(a) The restriction $H^{2}(X(\mu), \mathcal{F}) \rightarrow H^{2}(X(\lambda), \mathcal{F})$ has dense range;\\
(b) $H^i(X(\lambda), \mathcal{F})$ vanishes for all $i \geq 3$;\\
First we show that (a) holds. For this, we define $T \subseteq \mathrm{R}$ to be the set of all real numbers $\mu$ such that the restriction map
$$
H^{2}(X(\mu), \mathcal{F}) \rightarrow H^{2}(X(\lambda), \mathcal{F})
$$
has dense image for every real number $\lambda$ with $\lambda<\mu$.
Obviously, $T$ is not empty. In fact if $\mu_*:=\min \{f(y) ; y \in Y\}$, then clearly $]-\infty, \mu_*] \subset T$. To prove $T$ is open, we use the bumping method of Andreotti and Grauert. We fix some $\mu_0 \in T$. We shall find $\epsilon_o>0$ such that $\mu_o+\epsilon_o \in T$. For this, we consider Stein open subsets $U_i \subset \subset \mathbb{P}^{5}\setminus \nu(\mathbb{P}^{2}), i=1, \cdots, k$, such that $\left\{f=\mu_0\right\} \subset \bigcup_{i=1}^k U_i$ and choose functions $\left\{\theta_i\right\} \in C_0^{\infty}\left(U_i, \mathrm{R}\right), \theta_i \geq 0, i=1, \ldots, k$ with $\sum_{i=1}^k \theta_i(x)>0$ at any point $x \in\left\{\phi=\mu_0\right\}$. Define also smooth functions $f_j: \mathbb{P}^{5}\setminus \nu(\mathbb{P}^{2}) \rightarrow \mathrm{R}$ by

$$
f_j:=f-\sum_{i=1}^j c_i \theta_i, j=1, \ldots, k
$$

where $c_i>0$ are sufficiently small constants such that $f_o:=f, f_1, \ldots, f_k$, are $r$ convex with corners. Set

$$
X_j:=\left\{x \in  \mathbb{P}^{5}\setminus \nu(\mathbb{P}^{2}); f_j(x)<\mu_o\right\}, j=1, \ldots, k \text { and } X_o:=X(\mu_o) .
$$

Obviously, $X_j \backslash X_{j-1} \Subset U_j, X_j=X_{j-1} \cup(X_j \cap U_j)$ and $X_0 \subset \subset X_k$. Also since $f$ is proper, there exists $\varepsilon_o>0$ with $X(\mu_o+\epsilon_o) \subset X_k$. Furthermore, we remark that $dim_{\mathbb{C}}H^{3}(X_{j}, \mathcal{F})<\infty$
for all $j\in\{0,\cdots, k\}$. To see this, we consider the Mayer-Vietoris sequence for cohomology :
\begin{center}
$\cdots\rightarrow H^{3}(X_{j}, \mathcal{F})\rightarrow H^{3}(X_{j-1}, \mathcal{F})\oplus H^{3}(X_{j}\cap U_{j}, \mathcal{F})
\rightarrow H^{3}(X_{j-1}\cap U_{j}, \mathcal{F})\rightarrow\cdots$
\end{center}
Because $H^{3}(X_{j}\cap U_{j}, \mathcal{F})=H^{3}(X_{j-1}\cap U_{j}, \mathcal{F})=0$ by lemma $3$ for all $j=0,\cdots, k$, it follows that the restriction
$H^{3}(X_{k}, \mathcal{F})\rightarrow H^{3}(X_{0}, \mathcal{F})$ is surjective. Since in addition $X_{0}\subset\subset X_{k}$, we can conclude from
~\cite{ref2} that $dim_{\mathbb{C}}H^{3}(X_{j}, \mathcal{F})<\infty$ for $j=0,\cdots, k$.\\
We now consider the Mayer-Vietoris sequence for cohomology :
\begin{center}
$\cdots\rightarrow H^{2}(X_{j}, \mathcal{F})\rightarrow H^{2}(X_{j-1}, \mathcal{F})\oplus H^{2}(X_{j}\cap U_{j}, \mathcal{F})
\rightarrow H^{2}(X_{j-1}\cap U_{j}, \mathcal{F})\rightarrow H^{3}(X_{j}, \mathcal{F})\rightarrow H^{3}(X_{j-1}, \mathcal{F})\oplus H^{3}(X_{j}\cap U_{j}, \mathcal{F})\rightarrow H^{3}(X_{j-1}\cap U_{j}, \mathcal{F})\rightarrow\cdots$
\end{center}
It is easy to see that the restriction map $H^{3}(X_{j}, \mathcal{F})\rightarrow H^{3}(X_{j-1}, \mathcal{F})$
is an isomorphism. Therefore $H^{2}(X_{j}\cap U_{j}, \mathcal{F})\rightarrow H^{2}(X_{j-1}\cap U_{j}, \mathcal{F})$
is surjective, which implies according to the proof of Proposition $19$ in ~\cite{ref2} that the restriction map
$H^{2}(X_{j}, \mathcal{F})\rightarrow H^{2}(X_{j-1}, \mathcal{F})$ has dense range for $j=0,\cdots, k$.
It follows from the Mayer-Vietoris sequence for cohomology that the restriction map
$$
H^{2}\left(X\left(\mu_0+\varepsilon_0\right), \mathcal{F}\right) \rightarrow H^{2}(X(\mu), \mathcal{F})
$$
has dense image for all $\mu$ with $\mu_0 \leq \mu<\mu_0+\varepsilon_0$. Since $\mu_0 \in T$, then for every real number $\mu<\mu_0+\varepsilon_0$, the restriction
$$
H^{2}\left(X\left(\mu_0+\varepsilon_0\right), \mathcal{F}\right) \rightarrow H^{2}(X(\mu), \mathcal{F})
$$
has dense range, which shows that $\mu_0+\varepsilon_0 \in T$.
The set T is closed follows in a standard way from Proposition $20$ on page 246 in ~\cite{ref2}.
The proof of asertion (b) follows exactly the same steps as that of assertion (a), and will therefore be omitted.\\
\hspace*{.1in}In order to complete the proof of the theorem, note that for every
integer $j\geq 0$, we have $H^{r}(X(j), {\mathcal{F}})=0$ for all $r\geq 3$
and the restriction map
$$H^{r}(X(j+1), {\mathcal{F}})\rightarrow H^{r}(X(j), {\mathcal{F}})$$
has dense range if $r\geq 2$. Now the cohomological statement of theorem $3$
follows from (~\cite{ref2}, p. $250$).
\end{proof}
\section{$q$-convexity with corners}
\begin{Thm}
{Let $(n,q)$ be a pair of integers with $1\leq q\leq n-1$.
Then there exist an open subset $M\subset \mathbb{P}^{n}$ which
is $q$-complete with corners and a coherent analytic sheaf ${\mathcal{F}}$
on $\mathbb{P}^{n}$ such that $H^{\hat{q}}(M,{\mathcal{F}})\neq 0$, where
\[\hat{q}=n-[\frac{n-1}{q}]=\left\{\begin{array}{cc}
\tilde{q} \ \ if \ \  q\mid n\\
\tilde{q}-1 \ \ if \ \ q \nmid n
\end{array}
\right.\]}
Here $\tilde{q}=n-[\frac{n}{q}]+1$ and
$[\frac{n}{q}]$ is the integral part of $\frac{n}{q}$.
\end{Thm}
\begin{proof}
If $q$ divide $n$, it is easy to find $q$-complete with corners
complex manifolds which are not cohomologically $\hat{q}$-complete. (See e.g. ~\cite{ref11}).\\
Suppose now that $q\nmid n$, and consider the canonical quotient map
$$\begin{array}{ccccc}
\Pi & : & \mathbb{C}^{n+1}\setminus\{0\} & \to & \mathbb{P}^{n} \\
 & & z=(z_{0},\cdots,z_{n}) & \mapsto & \Pi(z)=[z] \\
\end{array}$$
Then clearly the sets\\
$A_{i}=\pi(\{z\in \mathbb{C}^{n+1}\setminus \{0\} : z_{iq}=\cdots z_{(i+1)q-1}=0\}) $ for 
$0\leq i\leq m-1$ and $A_{m}=\pi(\{z\in \mathbb{C}^{n+1}\setminus \{0\} : z_{mq}=\cdots z_{n}=0\}) $,
can be identified in a canonical way with the complex projective spaces $\mathbb{P}^{n-q}$ and $\mathbb{P}^{n-(r+1)}$, respectively,
where $n=mq+r$, with $m=[\frac{n}{q}]$ and $0<r<q$.
This implies that each $D_{i}=\mathbb{P}^{n}\setminus A_{i}$ is $q$-complete and for any $k$ with
$1\leq k\leq m,$ the set $D_{i_{1}}\cap\cdots\cap D_{i_{k}}$ is in particular $\hat{q}$-complete
for all $i_{1},\cdots, i_{k}\in\{0,1,\cdots, m\}$, since it is at worst $(mq-(m-1))$-complete and $mq-(m-1)\leq \hat{q}$.\\
\hspace*{.1in} On the other hand, the space $\mathbb{P}^{n}$ is not $n$-complete with corners, there exists a coherent analytic sheaf
$\mathcal{F}\in coh(\mathbb{P}^{n})$ such that $H^{n}(\mathbb{P}^{n}, \mathcal{F})\neq 0$.
Since $D_{i_{1}}\cap\cdots\cap D_{i_{k}}$ is $\hat{q}$-complete for all $k\in\{1,\cdots, m\}$ and $i_{1}, \cdots, i_{k}\in \{0,1,\cdots, m\}$, it follows from proposition $1$ of ~\cite{ref8} that
$$H^{\hat{q}}(D_{0}\cap\cdots\cap D_{m}, {\mathcal{F}})\cong
H^{n}(\mathbb{P}^{n}, {\mathcal{F}})\neq 0$$

\end{proof}
\noindent
$1$. Funding: Not applicable.\\
$2$. Informed Consent Statement: Not applicable.\\
$3$. Data Availability Statement: Not applicable.\\
$4$. Conflicts of Interest: The author declares no conflict of interest.\\
$5$. The name of all authors are written in full.

\vspace*{3pc}
\noindent
Corresponding author name's : Professor Youssef Alaoui\\
Hassan II Institute of Agronomy and Veterinary Sciences,\\
Madinat Al Irfane, BP 6202, Rabat, 10101, Morocco,\\
B.P.6202, Rabat-Instituts, 10101. Morocco.\\
Email : {\bf y.alaoui@iav.ac.ma or comp5123ster@gmail.com} \\\\
\end{document}